\documentclass[11pt]{article}
\usepackage{cite}
\usepackage{mathrsfs}
\usepackage{amsfonts}
\usepackage{amsmath}
\usepackage{amsfonts,amssymb,color}
\usepackage{dsfont}
\usepackage{curves}
\usepackage{mathrsfs}
\usepackage{pifont}
\usepackage{amssymb}
\usepackage{latexsym,amsmath,amssymb,amsfonts,epsfig,graphicx,cite,psfrag}
\usepackage{eepic,color,colordvi,amscd}

\newtheorem{theorem}{Theorem}[section]

\newtheorem{lemma}{Lemma}[section]
\newtheorem{corollary}{Corollary}[section]

\newtheorem{claim}{Claim}[section]
\newtheorem{conjecture}{Conjecture}[section]

\newtheorem{observation}{Observation}[section]

\newcommand{\qed}{\hfill\rule{0.5em}{0.809em}}

\def\emptyset{\mbox{{\rm \O}}}
\def\bar{\overline}

\def\qed{\hfill \rule{4pt}{7pt}}

\def\pf{\noindent {\it Proof. }}

\begin{document}

 \title{Coloring of some crown-free graphs}
  \author{Di Wu$^{1,}$\footnote{Email: 1975335772@qq.com},  \;\; Baogang  Xu$^{1,}$\footnote{Email: baogxu@njnu.edu.cn OR baogxu@hotmail.com. Supported by NSFC 11931006}\\\\
 	\small $^1$Institute of Mathematics, School of Mathematical Sciences\\
 	\small Nanjing Normal University, 1 Wenyuan Road,  Nanjing, 210023,  China}
 \date{}

 \maketitle
\begin{abstract}
Let $G$ and $H$ be two vertex disjoint graphs. The {\em union} $G\cup H$ is the graph with $V(G\cup H)=V(G)\cup (H)$ and $E(G\cup H)=E(G)\cup E(H)$. The {\em join} $G+H$ is the graph with $V(G+H)=V(G)+V(H)$ and $E(G+H)=E(G)\cup E(H)\cup\{xy\;|\; x\in V(G), y\in V(H)$$\}$. We use $P_k$ to denote a {\em path} on $k$ vertices, use {\em fork} to denote the graph obtained from $K_{1,3}$ by subdividing an edge once, and use {\em crown} to denote the graph $K_1+K_{1,3}$. In this paper, we show that (\romannumeral 1) $\chi(G)\le\frac{3}{2}(\omega^2(G)-\omega(G))$ if $G$ is (crown, $P_5$)-free, (\romannumeral 2) $\chi(G)\le\frac{1}{2}(\omega^2(G)+\omega(G))$ if $G$ is (crown, fork)-free, and (\romannumeral 3) $\chi(G)\le\frac{1}{2}\omega^2(G)+\frac{3}{2}\omega(G)+1$ if $G$ is (crown, $P_3\cup P_2$)-free.

\begin{flushleft}
\end{flushleft}

\end{abstract}

\newpage

\section{Introduction}
All graphs considered in this paper are finite and simple. We use $P_k$ and $C_k$ to denote a {\em path} and a {\em cycle} on $k$ vertices respectively, and follow \cite{BM08} for undefined notations and terminology. Let $G$ be a graph, and $X$ be a subset of $V(G)$. We use $G[X]$ to denote the subgraph of $G$ induced by $X$, and call $X$ a {\em clique} ({\em stable set}) if $G[X]$ is a complete graph (has no edge). The {\em clique number} $\omega(G)$ of $G$ is the maximum size taken over all cliques of $G$.

For $v\in V(G)$, let $N_G(v)$ be the set of vertices adjacent to $v$, $d_G(v)=|N_G(v)|$, $N_G[v]=N_G(v)\cup \{v\}$, $M_G(v)=V(G)\setminus N_G[v]$. For $X\subseteq V(G)$, let $N_G(X)=\{u\in V(G)\setminus X\;|\; u$ has a neighbor in $X\}$ and $M_G(X)=V(G)\setminus (X\cup N_G(X))$. If it does not cause any confusion, we will omit the subscript $G$ and simply write $N(v), d(v), N[v], M(v), N(X)$ and $M(X)$.  Let $\delta(G)$ denote the minimum degree of $G$.

Let $G$ and $H$ be two vertex disjoint graphs. The {\em union} $G\cup H$ is the graph with $V(G\cup H)=V(G)\cup (H)$ and $E(G\cup H)=E(G)\cup E(H)$. The {\em join} $G+H$ is the graph with $V(G+H)=V(G)+V(H)$ and $E(G+H)=E(G)\cup E(H)\cup\{xy\;|\; x\in V(G), y\in V(H)$$\}$. The union of $k$ copies of the same graph $G$ will be denoted by $kG$. The complement of a graph $G$ will be denoted by $\bar G$.  We say that $G$ induces $H$ if $G$ has an induced subgraph isomorphic to $H$, and say that $G$ is $H$-free otherwise. Analogously, for a family $\cal H$ of graphs, we say that $G$ is ${\cal H}$-free if $G$ induces no member of ${\cal H}$.

Let $G$ be a graph, $v\in V(G)$, and let $X$ and $Y$ be two subsets of $V(G)$. We say that $v$ is {\em complete} to $X$ if $v$ is adjacent to all vertices of $X$, and say that $v$ is {\em anticomplete} to $X$ if $v$ is not adjacent to any vertex of $X$. We say that $X$ is complete (resp. anticomplete) to $Y$ if each vertex of $X$ is complete (resp. anticomplete) to $Y$.

For positive integer $i$, let $N^i_G(X):=\{u\in V(G)\setminus X:$ min $\{d_G(u,v):v\in X\}=i\}$, where $d_G(u,v)$ is the distance between $u$ and $v$ in $G$. Then $N^1_G(X)=N_G(X)$ is the neighborhood of $X$ in $G$. Moreover, let $N^{\ge i}_G(X):= \cup^\infty_{j=i}N^i_G(X)$. We write $N^i_G(H)$ for $N^i_G(V(H))$. For integers $0\le i\le j$ and vertices $x\in N^i_G(H)$ and $y\in N^j_G(H)$, we call $x$ an {\em ancestor} of $y$ (resp. $y$ an {\em descendant of $x$}) if there exists a $uv$-path of length $j-i$.
When $G$ is clear from the context, we ignore the subscript $G$.

For $u, v\in V(G)$, we simply write $u\sim v$ if $uv\in E(G)$, and write $u\not\sim v$ if $uv\not\in E(G)$.  A {\em hole} of $G$ is an induced cycle of length at least 4, and a {\em $k$-hole} is a hole of length $k$. A $k$-hole is called an {\em odd hole} if $k$ is odd, and is called an {\em even hole} otherwise. An {\em antihole} is the complement of some hole. An odd (resp. even) antihole is defined analogously.

\begin{observation}\label{hole}
	The vertices of an odd hole or odd antihole cannot be divided into a stable set and a clique, or two cliques.
\end{observation}

Let $k$ be a positive integer, and let $[k]=\{1, 2, \ldots, k\}$. A $k$-{\em coloring} of $G$ is a mapping $c: V(G)\mapsto [k]$ such that $c(u)\neq c(v)$ whenever $u\sim v$ in $G$. The {\em chromatic number} $\chi(G)$ of $G$ is the minimum integer $k$ such that $G$ admits a $k$-coloring.  It is certain that $\chi(G)\ge \omega(G)$. A {\em perfect graph} is one such that  $\chi(H)=\omega(H)$ for all of its induced subgraphs $H$. The famous {\em Strong Perfect Graph Theorem}\cite{CRSR06} states that a graph is perfect if and only if it induces neither an odd hole nor an odd antihole. A graph is {\em perfectly divisible} \cite{HCT18} if for each of its induced subgraph $H$, $V(H)$ can be partitioned into $A$ and $B$ such that $H[A]$ is perfect and $\omega(H[B])<\omega(H)$. By a simple induction on $\omega(G)$ we have that $\chi(G)\le\frac{1}{2}(\omega^2(G)+\omega(G))$ for each perfectly divisible graph $G$.

\begin{figure}[htbp]\label{fig-1}
	\begin{center}
		\includegraphics[width=10cm]{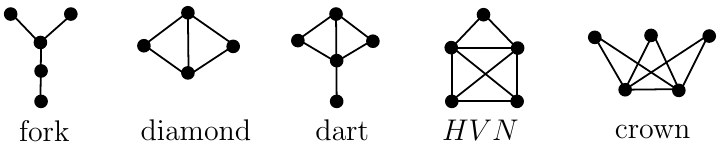}
	\end{center}
	\vskip -15pt
	\caption{Illustration of fork, diamond, dart, $HVN$, and crown.}
\end{figure}

A {\em fork} is the graph obtained from $K_{1,3}$ (usually called {\em claw}) by subdividing an edge once, a {\em diamond} is the graph $K_1+P_3$, a {\em dart} is the graph $K_1+(K_1\cup P_3)$, an $HVN$ is a $K_4$ together with one more
vertex which is adjacent to exactly two vertices of $K_4$, and a {\em crown} is the graph $K_1+K_{1,3}$, see Figure \ref{fig-1} for these configurations.

A family $\cal G$ of graphs is said to be $\chi$-{\em bounded} if there is a function $f$ such that $\chi(G)\le f(\omega(G))$ for every $G\in\cal G$, and if such a function does exist for $\cal G$, then $f$ is said to be a {\em binding function} of $\cal G$ \cite{Gy75}. Erd\"{o}s \cite{ER59} proved that for any positive integers $k,l\ge3$, there exists a graph $G$ with $\chi(G)\ge k$ and no cycles of length less than $l$. This result motivates us to study the chromatic number of $F$-free graphs, where $F$ is a forest (a disjoint union of trees). Gy\'{a}rf\'{a}s \cite{Gy75,Gy87} and Sumner \cite{Su81} independently, proposed the following conjecture.

\begin{conjecture}\label{tree}\cite{Gy75, Gy87, Su81}
	For every tree $T$, $T$-free graphs are $\chi$-bounded.
\end{conjecture}

Since $P_4$-free graphs are perfect, determining an optimal binding function of $P_5$-free graphs or fork-free graphs or $(P_3\cup P_2)$-free graphs attracts much attention. Choudum {\em et al}\cite{CKS07} (see also \cite{SR19,SS20}), and Karthick {\em et al}\cite{KKS22}, proposed the following two still open conjectures, respectively.

\begin{conjecture}\label{conj-2}\cite{CKS07}
	There exists a constant $c$ such that for every $P_5$-free graph $G$, $\chi(G)\le c\omega^2(G)$.
\end{conjecture}

\begin{conjecture}\label{conj-3}\cite{KKS22}
	The class of fork-free graphs is perfectly divisible.
\end{conjecture}

Up to now, the best known $\chi$-binding function for $P_5$-free graphs is $f(\omega)=\omega(G)^{\log_2\omega(G)}$ \cite{SS21}, the best known $\chi$-binding function for fork-free graphs is $f(\omega)=7.5\omega^2(G)$ \cite{LW21} and the best known $\chi$-binding function for $(P_3\cup P_2)$-free graphs is $f(\omega)=\frac{1}{6}\omega(G)(\omega(G)+1)(\omega(G)+2)$ \cite{BA18}. We refer the interested readers to \cite{GM22} for results of $P_5$-free graphs, and to \cite{RS04,SR19,SS20} for more results and problems about the $\chi$-bounded problem.

Recently, Dong {\em et al} \cite{DXX22} proved that $(P_5, C_5, K_{2,3})$-free graphs are perfectly divisible, and  $\chi(G)\le 2\omega^2(G)-\omega(G)-3$ if $G$ is $(P_5, K_{2,3})$-free with $\omega(G)\ge2$. Notice that crown can be obtained from $K_{2,3}$ by adding an edge joining two vertices of degree 3. Let $H\in\{P_5,$ fork, $P_3\cup P_2\}$. In this paper, we pay attention to (crown, $H$)-free graphs. Brause {\em  et al} \cite{BRS19} proved that the class of $(2K_2,3K_1)$-free graphs has no linear $\chi$-binding function and so one cannot expect a linear binding function for (crown, $H$)-free graphs.

Kim \cite{KIM95} showed that the Ramsey number $R(3,t)$ has order of magnitude $O(\frac{t^2}{\log t})$, and thus for any $K_{1,3}$-free graph $G$, $\chi(G)\le O(\frac{\omega^2(G)}{\log \omega(G)})$. Chudovsky and Seymour \cite{CS05} have given a detailed structural classification of $K_{1,3}$-free graphs, and they proved in \cite{CS10} that every connected $K_{1,3}$-free graph $G$ with a stable set of size at least three satisfies $\chi(G)\le2\omega(G)$. It is known from \cite{RB93} that for a fork-free graph $G$, $\chi(G)\le3$ if $\omega(G)=2$. Therefore, for a $K_{1,3}$-free graph $G$, $\chi(G)\le3=\frac{1}{2}(\omega^2(G)+\omega(G))$ if $\omega(G)=2$. Since $2\omega(G)\le\frac{1}{2}(\omega^2(G)+\omega(G))$ for any  $\omega(G)\ge3$, we have that every connected $K_{1,3}$-free graph $G$ with a stable set of size at least three satisfies $\chi(G)\le\frac{1}{2}(\omega^2(G)+\omega(G))$. So, since $3K_1$-free graphs are perfectly divisible, we have that

\begin{lemma}\label{claw}
	Every connected $K_{1,3}$-free graph $G$ satisfies $\chi(G)\le\frac{1}{2}(\omega^2(G)+\omega(G))$.
\end{lemma}

In this paper, we firstly present a relation on the binding functions of $(P_5,H)$-free graphs and $(P_5,K_1+H)$-free graphs, for some given graph $H$.

\begin{theorem}\label{P5}
	Let $H$ be a connected graph or the union of a connected graph and an isolated vertex with $|V(H)|\ge3$, and let $G$ be a connected $(P_5,K_1+H)$-free graph. If $(P_5,H)$-free graphs have a $\chi$-binding function $f(\omega)$, then $\chi(G)\le$ max$\{\chi(H),f(\omega(G)-1)\}+|V(H)|f(\omega(G)-1)$, and in particular, $\chi(G)\le (|V(H)|+1)f(\omega(G)-1)$ when $\chi(H)$ is sufficiently small.
\end{theorem}

As a consequence of Theorem \ref{P5}, we prove that

\begin{corollary}\label{P5*}
  $\chi(G)\le\frac{3}{2}(\omega^2(G)-\omega(G))$ if $G$ is (crown, $P_5$)-free.
\end{corollary}

Karthick {\em  et al} \cite{KKS22} proved that if $G$ is a connected (dart, fork)-free graph and contains an induced $K_{1,3}$ with the vertex $v$ of degree 3, then $G[M(v)]$ is perfect. In this paper, we prove a similar theorem on (crown, fork)-free graphs.

\begin{theorem}\label{Diamond}
	If $G$ is a connected (crown, fork)-free graph and contains an induced $K_{1,3}$ with the vertex $v$ of degree 3, then $G[M(v)]$ is perfect.
\end{theorem}

As its consequence, we prove that

\begin{corollary}\label{Diamond*}
	$\chi(G)\le\frac{1}{2}(\omega^2(G)+\omega(G))$ if $G$ is (crown, fork)-free.
\end{corollary}

For $(P_3\cup P_2)$-free graphs, the best known binding function is $f(\omega)=\frac{1}{6}\omega(G)(\omega(G)+1)(\omega(G)+2)$ \cite{BA18}. In \cite{PA22}, Prashant {\em  et al} proved that if $G$ is ($P_3\cup P_2$, diamond)-free, then $\chi(G)\le4$ when $\omega(G)=2$, $\chi(G)\le6$ when $\omega(G)=3$, $\chi(G)\le4$ when $\omega(G)=4$, and $G$ is perfect when $\omega(G)\ge5$, which improves some results of \cite{BA18}. They also proved \cite{PA22} that $\chi(G)\le\omega(G)+1$ if $G$ is a $(P_3\cup P_2, HVN)$-free graph with $\omega(G)\ge4$. In this paper, we prove that

\begin{theorem}\label{P3}
	$\chi(G)\le\frac{1}{2}\omega^2(G)+\frac{3}{2}\omega(G)+1$ if $G$ is (crown, $P_3\cup P_2$)-free.
\end{theorem}

\begin{figure}[htbp]\label{fig-2}
	\begin{center}
		\includegraphics[width=4cm]{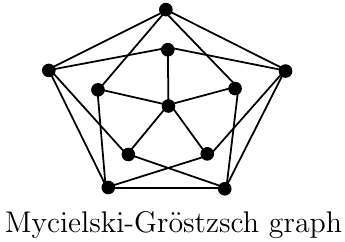}
	\end{center}
	\vskip -15pt
	\caption{Illustration of Mycielski-Gr\"{o}stzsch graph.}
\end{figure}

In \cite{WZ22}, Wang and Zhang proved that if $G$ is a $(P_3\cup P_2, K_3)$-free
graph, then $\chi(G)\le3$ unless $G$ is one of eight graphs with $\Delta(G)=5$ and $\chi(G)=4$, including the Mycielski-Gr\"{o}stzsch graph (see Figure \ref{fig-2}). It is certain that Mycielski-Gr\"{o}stzsch graph is also (crown, $P_3\cup P_2$)-free, and thus we cannot expect that (crown, $P_3\cup P_2$)-free graphs are perfectly divisible. In 1980, Wagon \cite{Wa80} proved that $\chi(G)\le\frac{1}{2}(\omega^2(G)+\omega(G))$ for every $2K_2$-free graph. It is reasonable to believe that the binding function of (crown, $P_3\cup P_2$)-free graphs can be improved to
$\frac{1}{2}\omega^2(G)+\frac{1}{2}\omega(G)+1$.

We will prove Theorem~\ref{P5} and Corollary~\ref{P5*} in Section 2, prove Theorem~\ref{Diamond} and Corollary~\ref{Diamond*} in Section 3, and prove Theorem~\ref{P3} in Section 4.

\section{(Crown, $P_5$)-free graphs }
Let $H$ be a connected graph, or the union of a connected graph and an isolated vertex. In this section, we consider $(P_5,K_1+H)$-free graphs with $|V(H)|\ge3$, and prove Theorem~\ref{P5} and Corollary~\ref{P5*}.

\medskip

\noindent{\em Proof of Theorem~\ref{P5}}: Let $f(\omega)$ be a binding function of $(P_5,H)$-free graphs, and let $G$ be a connected $(P_5,K_1+H)$-free graph. The theorem holds trivially if $G$ is $H$-free. So, we suppose that $G$ induces a $H$. Let $Q$ be an induced $H$ in $G$. It is certain that no vertex of $G$ is complete to $V(Q)$ as $G$ is $(K_1+H)$-free. We divide the proof process into two cases depending on $H$ is connected or not.

\medskip

\noindent{\bf Case} 1. Suppose that $H$ is a connected graph.

Since $|V(H)|\ge3$, it is certain that

\begin{equation}\label{eqa-1}
\mbox{for any vertex $u$ of $N(Q)$, $G[V(Q)\cup\{u\}]$ induces a $P_3$ starting from $u$}.
\end{equation}

Firstly, we prove that

\begin{equation}\label{eqa-2}
\mbox{$N^3(Q)=\emptyset$}.
\end{equation}

Suppose to its contrary, let $u_3\in N^3(Q)$. Then $G$ has an induced path $u_3u_2u_1u_0$, where $u_0\in V(Q)$ and $u_i\in N^i(Q)$ for $i\in \{1,2,3\}$. By (\ref{eqa-1}), $G[V(Q)\cup\{u_1\}]$ induces a $P_3$ starting from $u_1$, say $P$. Then $u_3u_2$-$P$ is an induced $P_5$ as $V(Q)$ is anticomplete to $\{u_2,u_3\}$, a contradiction. This proves (\ref{eqa-2}).

Let $u_1\in N(Q)$, and let $U$ be a component of $G[N^2(Q)]$ such that $u_1$ has a neighbor in $U$. If $u_1$ is not complete to $U$, then  $G[V(U)\cup\{u_1\}]$ induces a $P_3$ starting from $u_1$, say $P'$. By (\ref{eqa-1}), $G[V(Q)\cup\{u_1\}]$ induces a $P_3$ starting from $u_1$, say $P''$. Then $P'$-$u_1$-$P''$ is an induced $P_5$ as $V(Q)$ is anticomplete to $V(U)$, a contradiction. Therefore,

\begin{equation}\label{eqa-3}
\mbox{$u_1$ is complete to $V(U)$}.
\end{equation}

By (\ref{eqa-3}), we have that $U$ is $H$-free to forbid an induced $K_1+H$, and thus $\chi(U)\le f(\omega(G)-1)$. Consequently, $\chi(G[N^2(Q)])\le f(\omega(G)-1)$.

Now, for each vertex $u$ in $V(Q)$, $N(u)$ induces a $H$-free graph as otherwise $u$ and the $H$ induce a $K_1+H$. This proves that $\chi(G[N(Q)])\le |V(H)|f(\omega(G)-1)$.

Recall $N^3(Q)=\emptyset$, we have that $N^{\ge3}(Q)=\emptyset$. Since $\chi(G[N^2(Q)])\le f(\omega(G)-1)$ and $V(Q)$ is anticomplete to $N^2(Q)$, we have that $\chi(G)\le$ max$\{\chi(H),\chi(G[N^2(Q)])\}+\chi(G[N(Q)])\le$ max$\{\chi(H),f(\omega(G)-1)\}+|V(H)|f(\omega(G)-1)$.

\medskip

\noindent{\bf Case} 2. Suppose that $H$ is the union of a connected graph and an isolated vertex.

We still let $Q$ be an induced $H$ in $G$. Let $v'$ be the isolated vertex of $Q$, and $Q'=Q\setminus \{v'\}$. Since $G$ is connected, there exists a path from $v'$ to $Q'$. We choose the shortest one, say $P=v_1v_2\dots v_n$, where $v_1=v'$ and $v_n\in V(Q')$. We will prove that

\begin{equation}\label{c-1}
\mbox{$|V(P)|=3$}.
\end{equation}

It suffices to show $|V(P)|\le3$ as $v'$ is an isolated vertex in $H$. Suppose to its contrary that $|V(P)|\ge4$. Since $G$ is $P_5$-free, we may assume that $|V(P)|=4$, and let $P=v'v_2v_3v_4$. By the minimality of $P$, $\{v',v_2\}$ is anticomplete to $V(Q')$. If $v_3$ is complete to $V(Q')$, then $V(Q')\cup \{v_2,v_3\}$ induces a $K_1+H$, a contradiction. Therefore, $v_3$ is not complete to $V(Q')$, that is to say, $G[V(Q)\cup\{v_3\}]$ induces a $P_3$ starting from $v_3$, say $P'$. Now, $v'v_2v_3P'$ is an induced $P_5$, a contradiction. This proves (\ref{c-1}).

So, $P=v'v_2v_3$. It is certain that $v_2\in N(Q)$ and $v_3\in V(Q')$. Next, we prove that

\begin{equation}\label{c-2}
\mbox{each vertex of $N^2(Q)$ has an ancester in $V(Q')$ if $N^2(Q)\ne\emptyset$}.
\end{equation}

Suppose to its contrary that there exists a vertex $u$ of $N^2(Q)$ which has no ancester in $V(Q')$. Then the ancester of $u$ in $V(Q)$ is $v'$. Therefore, $G$ has an induced path $uu_1v'$, where $u_1\in N(Q)$. Since $v_2\sim v_3$ and $v_3$ cannot be an ancestor of $u$, we have that $u_1\ne v_2$. To forbid a $K_1+H$ induced by $V(Q')\cup \{v',v_2\}$, we have that $v_2$ is not complete to $V(Q')$. Thus $G[V(Q)\cup\{v_2\}]$ induces a $P_3$ starting from $v_2$, say $P''$. Since $u$ has no ancester in $V(Q')$, we have that $u\not\sim v_2$ and $u_1$ is anticomplete to $V(Q')$. Now, $uu_1v_2P''$ is an induced $P_5$ if $u_1\sim v_2$, and $uu_1v'v_2v_3$ is an induced $P_5$ if $u_1\not\sim v_2$, both are contradictions. This proves (\ref{c-2}).

Let $u'\in N(Q)$, suppose that $u'$ has a neighbor, say $u_2'$ in $N^2(Q)$. If $u'$ is complete to $V(Q')$, then $V(Q')\cup \{u',u'_2\}$ induces a $K_1+H$, which leads to a contradiction. So,

\begin{equation}\label{c-3}
\mbox{$u'$ is not complete to $V(Q')$ for each $u'\in N(Q)$ which has a neighbor in $N^2(Q)$}.
\end{equation}

Since $G$ is $P_5$-free, we can deduce that $N^{\ge3}(Q)=\emptyset$. Let $T$ be a component of $G[N^2(Q)]$, and let $t$ be a neighbor of $T$ in $N(Q)$. Then $t$ is complete to $T$ as otherwise one can find an induced $P_5$ by (\ref{c-2}) and (\ref{c-3}). Consequently, we have that $\chi(T)\le f(\omega(G)-1)$, and it follows that $\chi(G[N^2(Q)])\le f(\omega(G)-1)$. Therefore, $\chi(G)\le$ max$\{\chi(H),\chi(G[N^2(Q)])\}+\chi(G[N(Q)])\le$ max$\{\chi(H),f(\omega(G)-1)\}+|V(H)|f(\omega(G)-1)$.

This completes the proof of Theorem \ref{P5}.
\qed

\medskip

Let $G$ be a ($P_5$, crown)-free graph. Note that $(P_5,K_{1,3})$-free graphs have a $\chi$-binding function $f(\omega)=\frac{1}{2}(\omega^2+\omega)$ by Lemma \ref{claw}. We may assume that $G$ contains an induced $K_{1,3}$, say $Q$, and $V(Q)=\{v_1,v_2,v_3,v_4\}$, where $d_Q(v_1)=3$. By Theorem \ref{P5}, we may color $G$ by $\frac{5}{2}(\omega^2(G)-\omega(G))$ colors as $Q$ is a connected graph and $|V(Q)|=4$. We may improve this bound to $\frac{3}{2}(\omega^2(G)-\omega(G))$.

\medskip

\noindent{\em Proof of  Corollary~\ref{P5*}}: Without loss of generality, we suppose that $G$ induces a $K_{1,3}$, say $Q$. If $N^3(Q)\ne\emptyset$, then one can find an induced $P_5$ or an induced crown, a contradiction. So, $N^3(Q)=\emptyset$.
We only need to color $G[V(Q)\cup N(Q)\cup N^2(Q)]$ as $G=G[V(Q)\cup N(Q)\cup N^2(Q)]$. We prove firstly that

\begin{equation}\label{co-1}
\mbox{if $(N(v_i)\setminus \{v_1\})\subseteq N(V(Q)\setminus\{v_i\})$ for some $i\in\{2,3,4\}$, then $\chi(G)\le\frac{3}{2}(\omega^2(G)-\omega(G))$}.
\end{equation}

By symmetry, we may assume that $N(v_4)\setminus\{v_1\}\subseteq N(\{v_1,v_2,v_3\})$. Since $G$ is crown-free, we have that $N(v)$ is $K_{1,3}$-free for each $v\in V(G)$, which implies that $\chi(N(v))\le \frac{1}{2}(\omega^2(G)-\omega(G))$ for each $v\in V(G)$ by Lemma \ref{claw}.

Firstly, we color $N(v_1)$ and $N(v_2)$ by $(\omega^2(G)-\omega(G))$ colors. Consequently, by (\ref{eqa-3}), we can color $N^2(Q)$ by $\frac{1}{2}(\omega^2(G)-\omega(G))$ other colors. Then, it suffices to color $V(G)\setminus (N(v_1)\cup N(v_2)\cup N^2(Q))$. It is certain that  $V(G)\setminus (N(v_1)\cup N(v_2)\cup N^2(Q))\subseteq N(Q)$ as $N^3(Q)=\emptyset$.

Let $N_3=N(v_3)\setminus (N(v_1)\cup N(v_2))$. Then $N_3\subseteq N(Q)$ and $N_3$ is anticomplete to $N^2(Q)$ as otherwise there exists $u\in N_3$ and $u_2\in N^2(Q)$ such that $u\sim u_2$, which implies that $\{u_2,u,v_1,v_2,v_3\}$ induces a $P_5$, a contradiction. Therefore, we may color $N_3$ by the $\frac{1}{2}(\omega^2(G)-\omega(G))$ colors of $N^2(Q)$. Since $V(G)\setminus (N(v_1)\cup N(v_2)\cup N^2(Q)\cup N_3)\subseteq (N(v_4)\setminus \{v_1\})\subseteq N(\{v_1,v_2,v_3\})\subseteq (N(v_1)\cup N(v_2)\cup N_3)$, we have colored $V(G)\setminus (N(v_1)\cup N(v_2)\cup N^2(Q))$ by the $\frac{3}{2}(\omega^2(G)-\omega(G))$ colors of $N(v_1)\cup N(v_2)\cup N_3$. Therefore, we may color $V(Q)\cup N(Q)\cup N^2(Q)$ by $\frac{3}{2}(\omega^2(G)-\omega(G))$ colors, which implies that $\chi(G)\le\frac{3}{2}(\omega^2(G)-\omega(G))$. This proves (\ref{co-1}).

Next, we may assume by symmetry that there exists a vertex $v_2'\in N(v_2)\setminus \{v_1\}$ such that $v_2'\notin N(\{v_1,v_3,v_4\})$. We still define that
$N_3=N(v_3)\setminus (N(v_1)\cup N(v_2))$. By (\ref{co-1}), $N_3\ne\emptyset$. Let $v_3'\in N_3$. Since $v_2'$ is anticomplete to $\{v_1,v_3\}$ and $v_3'$ is anticomplete to $\{v_1,v_2\}$, we have that $v_2'\sim v_3'$ as otherwise $\{v_3',v_1,v_2,v_3,v_2'\}$ induces a $P_5$. But now, $v_3'\sim v_4$ to forbid an induced $P_5$ on $\{v_2',v_3',v_1,v_3,v_4\}$ as $v_2'$ is anticomplete to $\{v_1,v_3,v_4\}$ and $v_3'\sim v_1$. Therefore $v_3'\in N(v_4)\setminus \{v_1\}$. So, $N_3\subseteq N(v_4)\setminus \{v_1\}$, which implies that $(N(v_3)\setminus \{v_1\})\subseteq N(\{v_1,v_2,v_4\})$. By (\ref{co-1}), $\chi(G)\le\frac{3}{2}(\omega^2(G)-\omega(G))$. This proves Corollary \ref{P5*}.
\qed

\section{(Crown, fork)-free graphs }

In this section, we turn our attention to (crown, fork)-free graphs, and prove Theorem~\ref{Diamond} and Corollary~\ref{Diamond*}. Let $G$ be a  fork-free graph, and $C=v_1v_2\dots v_nv_1$ an odd hole contained in $G$, and $uv$ an edge of $G- V(C)$. If $v$ is anticomplete to $V(C)$ and $N(u)\cap V(C)=\{v_j, v_{j+1}\}$ for some $j$, then we call $v_jv_{j+1}$ a {\em bad edge} of $u$.

The following two lemmas will be used in our proof.

\begin{lemma}\label{le-1} ({Lemma 2.4 of \cite{WX22}})
	Let $G$ be a  fork-free graph, and $C=v_1v_2\dots v_nv_1$ an odd hole contained in $G$. If $G- V(C)$ has an edge $uv$ such that $u$ is not anticomplete to $V(C)$ but $v$ is, then either $u$ is complete to $V(C)$, or $v_jv_{j+1}$, for some $j\in\{1,2,\dots,n\}$ is a bad edge of $u$.
\end{lemma}

\begin{lemma}\label{le-2}
	Let $G$ be a fork-free graph, and $C$ an odd antihole  with $|V(C)|\ge7$ contained in $G$. If $G- V(C)$ has an edge $uv$ such that $u$ is not anticomplete to $V(C)$ but $v$ is, then $u$ is complete to $V(C)$.
\end{lemma}
\pf Let $C=\{v_1, v_2,\ldots, v_n\}$ such that $n\ge7$ and $v_i\sim v_j$ whenever $|i-j|\ne 1$ (indices are modulo
$n$). Let $uv$ be an edge of $G- V(C)$ such that $u$ is not anticomplete to $V(C)$ but $v$ is, and suppose that $u$ is not complete to $V(C)$.

If $u$ has two consecutive neighbors in $C$, then there is some $i$ such that $v_i, v_{i+1}\in N(u)$ and $v_{i+2}\notin N(u)$, which implies a fork induced by $\{v_{i+2}, v_i, u, v_{i+1}, v\}$. This shows that $u$ does not have two consecutive neighbors in $C$.

Without loss of generality, suppose that $u\sim v_1$. Then $u\not\sim v_2$ and $u\not\sim v_n$.

If $u\sim v_3$, then $u\not\sim v_4$, and thus $u\sim v_5$ to avoid an induced fork on $\{v,u,v_1,v_4,v_5\}$. Consequently $u\sim v_j$ for $j\in\{3,5,7,\dots,n-2\}$, and $u\not\sim v_k$ for $k\in\{2,4,6,\dots,n-1\}$. In particular, $u\sim v_{n-4}$ and $u\not\sim v_{n-1}$. If $u\sim v_n$, then $u$ has two consecutive neighbors $v_1$ and $v_n$. If $u\not\sim v_n$, then $\{v,u,v_{n-4},v_{n-1},v_n\}$ induces a fork. Both are contradictions. Therefore, $u\not\sim v_3$. To forbid a fork on $\{a,b,v_1,v_3,v_4\}$, we have that $u\sim v_4$, and thus $u\not\sim v_5$ as $u$ has no two consecutive neighbors in $C$. If $u\not\sim v_6$  then $\{v,u,v_1,v_5,v_6\}$ induces a fork. If $u\sim v_6$, then $\{v,u,v_6,v_2,v_3\}$ induces a fork. Both are contradictions. This proves Lemma \ref{le-2}.
\qed

\medskip

\noindent{\em Proof of  Theorem~\ref{Diamond}}: Let $G$ be a (crown, fork)-free graph. By Lemma \ref{claw}, we may assume that $G$ contains an induced $K_{1,3}$, say $Q$. Let $V(Q)=\{v_1,v_2,v_3,v_4\}$, where $d_Q(v_1)=3$. It is certain that no vertex of $V(G)\setminus V(Q)$ is complete to $V(Q)$ as $G$ is crown-free. We prove firstly that

\begin{equation}\label{t2-1}
\mbox{$N^3(Q)=\emptyset$}.
\end{equation}

Suppose to its contrary that there exists a vertex $u_3$ in $N^3(Q)$. Then $G$ has an induced path $u_3u_2u_1u_0$, where $u_i\in N^i(Q)$ for $i\in\{1,2,3\}$ and $u_0\in V(Q)$.

Suppose $u_1\sim v_1$. If $u_1$ is adjacent to at most one element of $\{v_2,v_3,v_4\}$, we may assume that $u_1\not\sim v_3$ and $u_1\not\sim v_4$, then $\{u_1,u_2,v_1,v_3,v_4\}$ induces a fork, a contradiction. Otherwise, $u_1$ is adjacent to at least two elements of $\{v_2,v_3,v_4\}$, we may assume that $u_1\sim v_3$ and $u_1\sim v_4$, which implies that $\{u_1,u_2,u_3,v_3,v_4\}$ induces a fork, a contradiction. Therefore, $u_1\not\sim v_1$. By symmetry, we suppose $u_1\sim v_2$.

If $u_1$ is anticomplete to $\{v_3,v_4\}$, then $\{u_1,v_1,v_2,v_3,v_4\}$ induces a fork, which leads to a contradiction. Therefore, $u_1$ is adjacent to at least one element of $\{v_3,v_4\}$. Without loss of generality, suppose that $u_1\sim v_3$. Then $\{u_1,u_2,u_3,v_2,v_3\}$ induces a fork, a contradiction. This proves (\ref{t2-1}).

Next, we prove that

\begin{equation}\label{t2-2}
\mbox{each component of $N^2(Q)$ is a clique}.
\end{equation}

Suppose that $N^2(Q)\ne\emptyset$ and $U$ is a component in $N^2(Q)$. It is certain that $U$ is anticomplete to $V(Q)\cup (N^2(Q)\setminus V(U))$.

Suppose that there exists a vertex $t_1$ in $N(Q)$ such that $t_1$ is neither complete to $U$ nor anticomplete to $U$. Then $U$ has two adjacent vertices, say $t_2$ and $t_2'$, such that $t_1\sim t_2$ and $t_1\not\sim t_2'$. By relacing $t_2',t_2,t_1,t_0$ with $u_3,u_2,u_1,u_0$, respectively, we can find an induced a fork with the same arguments as that used in proving (\ref{t2-1}), which leads to a contradiction. Therefore, we have that each vertex of $N(Q)$ is either complete to $U$ or anticomplete to $U$.

Suppose that $U$ is not a clique. Let $t_3$ and $t_3'$ be two nonadjacent vertices of $U$, and let $s$ be a neighbor of $U$ in $N(Q)$. It is certain that $s\sim t_3$ and $s\sim t_3'$. Since $s$ is not complete to $V(Q)$, there exists an induced $P_2$ starting from $s$ in $G[V(Q)\cup\{s\}]$, say $P$. Then $\{t_3,t_3'\}\cup V(P)$ induces a fork, a contradiction. So, $U$ is a clique. This proves (\ref{t2-2}).

\medskip

Let $N_1=N(v_1)\cap N(Q), N_2=\{v\in N(Q)\setminus N_1 | N(v)\cap V(Q)=\{v_2,v_3 \}\},N_3=\{v\in N(Q)\setminus N_1 | N(v)\cap V(Q)=\{v_2,v_4 \}\}, N_4=\{v\in N(Q)\setminus N_1 | N(v)\cap V(Q)=\{v_3,v_4 \}\}$, and $N_5=\{v\in N(Q)\setminus N_1 | N(v)\cap V(Q)=\{v_2,v_3,v_4 \}\}$. Since $G$ is fork-free, we have that no vertex of $N(Q)\setminus N_1$ can have a single neighbor in $V(Q)$, and thus $N(Q)=\cup^5_{i=1}N_i$.

Since $G$ is  crown-free, we have that $N_5$ is a stable set. If $N_2$ has two nonadjacent vertices, say $u$ and $v$, then $\{u,v,v_1,v_2,v_4\}$ induces a fork. The same contradiction happens if $N_3$ or $N_4$ has two nonadjacent vertices. Therefore,

\begin{equation}\label{t2-3}
\mbox{$N_i$ is a clique for each $i\in\{2,3,4\}$, and $N_5$ is a stable set}.
\end{equation}

To complete the proof of Theorem~\ref{Diamond}, we will show $G[M(v_1)]$ is perfect.

\begin{claim}\label{t2-4}
\mbox{$G[(N(Q)\setminus N_1)\cup N^2(Q)]$ contains no odd hole}.
\end{claim}
\pf Suppose to its contrary that $C=x_1x_2\dots x_nx_1$ is an odd hole contained in $G[(N(Q)\setminus N_1)\cup N^2(Q)]$. Without loss of generality, we suppose that $v_2$ has a neighbor in $C$.  Notice that $v_1$ is anticomplete to $V(C)$.  By Lemma \ref{le-1}, we have that either $N(v_2)\cap V(C)= V(C)$, or $N(v_2)\cap V(C)=\{x_j, x_{j+1}\}$ for some $j\in \{1,2,\dots,n\}$.

Suppose that $V(C)\cap N^2(Q)\ne\emptyset$. Then $v_2$ has a bad edge, say $e=u'v'$, in $C$ since $v_2$ is anticomplete to $N^2(Q)$. It is certain that $\{u',v'\}\subseteq N_2\cup N_3\cup N_5$ as $v_2$ is anticomplete to $N_4$. Moreover, $V(C)\cap ((N_2\cup N_3\cup N_5)\setminus\{u',v'\})=\emptyset$ as $v_2$ is complete to $N_2\cup N_3\cup N_5$.

If $\{u',v'\}\subseteq N_2\cup N_5$, then $e$ is also a bad edge of $v_3$ in $C$ as $v_3$ is complete to $N_2\cup N_5$. Since $v_3$ is complete to $N_4$, it follows that $V(C)\cap N_4=\emptyset$. Therefore, $V(C)$ is contained in $\{u',v'\}\cup N^2(Q)$ as $V(C)\cap ((N_2\cup N_3\cup N_5)\setminus\{u',v'\})=\emptyset$, contradicting (\ref{t2-2}). So, $\{u',v'\}\not\subseteq N_2\cup N_5$. By symmetry, we have $\{u',v'\}\not\subseteq N_3\cup N_5$ also, and thus
$$\{u',v'\}\subseteq N_2\cup N_3 \mbox{ and } |\{u',v'\}\cap N_2|=|\{u',v'\}\cap N_3|=1.$$

By symmetry, we may assume that $u'\in N_2$ and $v'\in N_3$. If there exists a vertex $w\in N(C)\cap N_4$, then $u'w$ is a bad edge of $v_3$ as $v_3\sim u'$ and $v_3\sim w$, and similarly $v'w$ is a bad edge of $v_4$ as $v_4\sim v'$ and $ v_4\sim w$, which leads to a contradiction that $\{u',v',w\}\subseteq V(C)$ but induces a triangle. Therefore, $V(C)\cap N_4=\emptyset$. But then, $V(C)\subseteq \{u',v'\}\cup N^2(Q)$, contradicting (\ref{t2-2}) again. This proves that $V(C)\cap N^2(Q)=\emptyset$ and $V(C)\subseteq N(Q)\setminus N_1=\cup^5_{i=2}N_i$.

Recall that an odd hole cannot be divided into two cliques by Observation \ref{hole}. By (\ref{t2-3}), we have that either $V(C)\cap N_4\ne\emptyset$ or $V(C)\cap N_5\ne\emptyset$ as otherwise $V(C)$ is divided into two cliques $G[N_2]$ and $G[N_3]$, which leads to a contradiction. By the definition of $N_4$ and $N_5$, we have that $v_3$ and $v_4$ both have a neighbor in $V(C)$.

Since $N_4$ is a clique, it follows that $|V(C)\cap N_4|\le2$. So, $|V(C)\cap (N_2\cup N_3\cup N_5)|\ge3$ as $n\ge5$. Moreover, by Lemma \ref{le-1}, $v_2$ is complete to $V(C)$ as $v_2$ is complete to $N_2\cup N_3\cup N_5$. By symmetry,  $\{v_3, v_4\}$ is complete to $V(C)$ as $v_3$ and $v_4$ both have a neighbor in $V(C)$.

But now, we have that $V(C)\cap N_2=\emptyset$ as $v_4$ is anticomplete to $N_2$, $V(C)\cap N_3=\emptyset$ as $v_3$ is anticomplete to $N_3$, and $V(C)\cap N_4=\emptyset$ as $v_2$ is anticomplete to $N_4$. So, $V(C)\subseteq N_5$, contradicting (\ref{t2-3}). This proves Claim~\ref{t2-4}. \qed

\begin{claim}\label{t2-5}
\mbox{$(N(Q)\setminus N_1)\cup N^2(Q)$ contains no odd antihole}.
\end{claim}
\pf
Suppose to its contrary that $C'$ is an odd antihole with $V(C')$ contained in $(N(Q)\setminus N_1)\cup N^2(Q)$. Without loss of generality, we suppose that $v_2$ has a neighbor in $C'$.  Notice that $v_1$ is anticomplete to $V(C')$.  By Lemma \ref{le-2}, we have that $N(v_2)\cap V(C')= V(C')$.

Recall that an odd antihole cannot be divided into two cliques by Observation \ref{hole}. By (\ref{t2-3}), we have that either $V(C')\cap N_4\ne\emptyset$ or $V(C')\cap N_5\ne\emptyset$ as otherwise $V(C')$ is divided into two cliques $G[N_2]$ and $G[N_3]$, which leads to a contradiction. By the definition of $N_4$ and $N_5$, we have that  $v_3$ and $v_4$ both have a neighbor in $V(C')$. It follows that $\{v_3, v_4\}$ is complete to $V(C')$ by Lemma \ref{le-2}.

But now, we have that $V(C')\cap N_2=\emptyset$ as $v_4$ is anticomplete to $N_2$, $V(C')\cap N_3=\emptyset$ as $v_3$ is anticomplete to $N_3$, and $V(C')\cap N_4=\emptyset$ as $v_2$ is anticomplete to $N_4$. So, $V(C')\subseteq N_5$, contradicting (\ref{t2-3}).  This proves Claim~\ref{t2-5}. \qed

\medskip

By Claims~\ref{t2-4} and \ref{t2-5}, and by the Strong Perfect Graph Theorem, we have that $(N(Q)\setminus N_1)\cup N^2(Q)$ is perfect. Note that $N(v_1)=N_1\cup\{v_2,v_3,v_4\}$ and $M(v_1)=(N(Q)\setminus N_1)\cup N^2(Q)$. This completes the proof of Theorem~\ref{Diamond}. \qed

\medskip

\noindent{\em Proof of  Corollary~\ref{Diamond*}}: Suppose that the statement does not hold. Let $G$ be a minimal (crown, fork)-free graph with $\chi(G)>\frac{1}{2}(\omega^2(G)+\omega(G))$. By Lemma \ref{claw}, we have that $G$ must contain an induced $K_{1,3}$. Let $v$ be the vertex of degree 3 in an induced $K_{1,3}$ of $G$. By Theorem~\ref{Diamond}, $M(v)$ is perfect. Since $\chi(N(v))\le\frac{1}{2}((\omega(G)-1)^2+\omega(G)-1)$, we have that $\chi(G)\le\frac{1}{2}((\omega(G)-1)^2+\omega(G)-1)+\omega(G)=\frac{1}{2}(\omega^2(G)+\omega(G))$, a contradiction. Therefore, $\chi(G)\le\frac{1}{2}(\omega^2(G)+\omega(G))$ for each ($K_1+K_{1,3}$, fork)-free graph. This proves Corollary \ref{Diamond*}.\qed

\section{(Crown, $P_3\cup P_2$)-free graphs }

In this section, we focus on (crown, $P_3\cup P_2$)-free graphs, and prove Theorem \ref{P3}. Up to now, the best known binding function for $(P_3\cup P_2)$-free graph is due to Bharathi and Choudum.

\begin{theorem}\label{P3P2}\cite{BA18}
	If $G$ is $(P_3\cup P_2)$-free, then $\chi(G)\le\frac{1}{6}\omega(G)(\omega(G)+1)(\omega(G)+2)$.
\end{theorem}

Consequently, $\chi(G)\le4$ if $\omega(G)=2$ and $\chi(G)\le10$ if $\omega(G)=3$. This means that Theorem \ref{P3} holds when $\omega(G)\le3$. Therefore, we may assume that $\omega(G)\ge4$.

\medskip

\noindent{\em Proof of  Theorem~\ref{P3}}: Let $G$ be a (Crown, $P_3\cup P_2$)-free graph, and suppose that $\omega(G)\ge4$. By Lemma \ref{claw}, we may assume that $G$ contains an induced $K_{1,3}$, say $Q$. Let $V(Q)=\{v_1,v_2,v_3,v_4\}$, where $d_Q(v_1)=3$. It is certain that no vertex of $V(G)\setminus V(Q)$ is complete to $V(Q)$ as $G$ is crown-free. Since $G$ is $(P_3\cup P_2)$-free, we have that $N^{\ge2}(Q)$ is a stable set.

Let $N_1=N(v_1), N_2=\{v\in N(Q)\setminus N_1 | N(v)\cap V(Q)=\{v_2,v_3 \}\},N_3=\{v\in N(Q)\setminus N_1 | N(v)\cap V(Q)=\{v_2,v_4 \}\}, N_4=\{v\in N(Q)\setminus N_1 | N(v)\cap V(Q)=\{v_3,v_4 \}\}$, and $N_5=\{v\in N(Q)\setminus N_1 | N(v)\cap V(Q)=\{v_2,v_3,v_4 \}\}, N_6=\{v\in N(Q)\setminus N_1 | N(v)\cap V(Q)=\{v_4\}\}, N_7=\{v\in N(Q)\setminus N_1 | N(v)\cap V(Q)=\{v_3\}\}, N_8=\{v\in N(Q)\setminus N_1 | N(v)\cap V(Q)=\{v_2\}\}$. It is certain that $N(Q)\cup\{v_2,v_3,v_4\}=\cup^8_{i=1}N_i$.

Since $\{v_1,v_4\}$ is anticomplete to $N_2$, we have that $G[N_2]$ is $P_3$-free as $G$ is $(P_3\cup P_2)$-free. Then $N_2$ is a disjoint union of cliques. Let $x$ and $y$ be two adjacent vertices of $N_6$. Then $\{x,y,v_1,v_2,v_3\}$ induces a $P_3\cup P_2$, a contradiction. Therefore, $N_6$ is a stable set. By symmetry, $N_i$ is a disjoint union of cliques and $N_j$ is a stable set, for $i\in\{3,4\}$ and $j\in\{7,8\}$. Moreover, to forbid a crown, we have that $N_5$ is a stable set.

Let $u$ be a vertex in $N_2$, and $T$ be a clique of $G[N_3]$ with $|V(T)|\ge2$. If $T$ has an edge $u_1u_2$ such that $u$ is anticomplete to $\{u_1,u_2\}$, then $\{u,v_3,v_1,u_1,u_2\}$ induces a $P_3\cup P_2$. This contradiction implies that $|N(u)\cap V(T)|\ge |V(T)|-1$, and by symmetry, this inequality still holds when $T$ is a clique of $G[N_4]$. Therefore, let $u$ be a vertex of $N_i$ and $T$ be a clique of $G[N_j]$ with $|V(T)|\ge2$, for $i\in\{2,3,4\}$ and $j\in\{2,3,4\}\setminus\{i\}$, we may assume that

\begin{equation}\label{t3-1}
\mbox{$|N(u)\cap V(T)|\ge |V(T)|-1$}.
\end{equation}
We need the following two claims to complete our proof.

\begin{claim}\label{t3-2}
Let $i,j \in\{2,3,4\}$ and $i\ne j$. Then the following statements hold.

(1) If $\omega(G[N_i])\ge2$, then $N_{i+4}=\emptyset$.

\medskip

(2) If $\omega(G[N_i])\ge4$, then $N_j$ is a clique.

\medskip

(3) If $\omega(G[N_i])\ge3$ and $\omega(G[N_j])\ge3$, then both $N_i$ and $N_j$ are cliques.

\end{claim}

\pf We may by symmetry set $i=2$ and $j=3$. Let $T$ be a clique in $G[N_2]$.

Suppose $|V(T)|\ge2$. Let $u$ and $v$ be two adjacent vertices in $T$, and let $w$ be a vertex in $N_6$. To forbid a $P_3\cup P_2$ on $\{v_2,v_3,u,v_4,w\}$ or $\{v_2,v_3,v,v_4,w\}$, we have that $w$ is complete to $\{u,v\}$ as $w$ is anticomplete to $\{v_2,v_3\}$ and $v_4$ is anticomplete to $\{v_2,v_3,u,v\}$. But now, $\{v_2,v_3,u,v,w\}$ induces a crown, which leads a contradiction. So (1) holds.

Suppose $|V(T)|\ge4$ and $G[N_3]$ has two nonadjacent vertices $u_1$ and $u_2$. By (\ref{t3-1}), $|N(u_t)\cap V(T)|\ge |V(T)|-1\ge3$, for $t\in\{1,2\}$. Then $T$ must have an edge $u'v'$ such that $\{u_1,u_2\}$ is complete to $\{u',v'\}$. But now, $\{v_3,u_1,u_2,u',v'\}$ induces a crown as $\{v_3,u',v'\}$ is a stable set, which leads to a contradiction. So (2) holds.

Let $T$ be a clique with $|V(T)|=3$ contained in a component $T'$ of $G[N_2]$, and $S$ be a clique with $|V(S)|=3$ contained in a component $S'$ of $G[N_3]$. Let $V(T)=\{u_1',u_2',u_3'\}$ and $V(S)=\{v_1',v_2',v_3'\}$.

Suppose $u_4'\in N_2\setminus V(T')\ne\emptyset$. By (\ref{t3-1}), $u_4'$ must be complete to some two vertices of $V(S)$, by symmetry, $u_4'$ is complete to $\{v_1',v_2'\}$. By (\ref{t3-1}) again, $v_1'$ must be complete to some two vertices of $V(T)$ and $v_2'$ must be complete to some two vertices of $V(T)$. Therefore, $V(T)$ has at least one vertex which is complete to $\{v_1',v_2'\}$ as $|V(T)|=3$, by symmetry, $u_1'$ is complete to $\{v_1',v_2'\}$. But now, $\{v_1',v_2',v_4,u_1',u_4'\}$ induces a crown. This contradiction implies that $N_2\setminus V(T')=\emptyset$. By symmetry, $N_3\setminus V(S')=\emptyset$. So (3) holds.\qed

\begin{claim}\label{t3-3}
	
$\chi(G[\cup^8_{i=2}N_i])\le2\omega(G)$.
	
\end{claim}
\pf By Claim \ref{t3-2} (1), we have that $\chi(G[N_i\cup N_{i+4}])\le2$ if $\omega(G[N_i])\le2$, for $i\in\{2,3,4\}$. Without loss of generality, we suppose $\omega(G[N_2])\ge\omega(G[N_3])\ge\omega(G[N_4])$. Now, we divide the proof process into two cases depending on $\omega(G[N_2])\ge4$ or $\omega(G[N_2])\le3$.

\medskip

\noindent{\bf Case} 1. $\omega(G[N_2])\ge4$.

By the definition of $N_2$ and Claim \ref{t3-2} (2), we have that $\omega(G)\ge5$ and $N_j$ is a clique, for $j\in\{3,4\}$.

Suppose $\omega(G[N_3])\ge3$. By Claim \ref{t3-2} (3), $N_2$ is also a clique. Notice that an odd hole or odd antihole cannot be divided into two cliques by Observation \ref{hole}. By the Strong Perfect Graph Theorem, we have that $G[N_2\cup N_3]$ is perfect, which implies that $\chi(G[N_2\cup N_3])\le\omega(G)$. Moreover, by Claim \ref{t3-2} (1), $N_6=N_7=\emptyset$.

If $\omega(G[N_4])=1$, then $\chi(G[N_4\cup N_5\cup N_8])\le3$ as $N_5$ and $N_8$ are both stable sets. If $\omega(G[N_4])\ge2$, then $N_8=\emptyset$ by Claim \ref{t3-2} (1), and $G[N_4\cup N_5]$ is perfect by Observation \ref{hole}, which implies that $\chi(G[N_4\cup N_5\cup N_8])\le\omega(G)-1$. These follow that $\chi(G[\cup^8_{i=2}N_i])\le\chi(G[N_2\cup N_3])+\chi(G[N_4\cup N_5\cup N_8])\le2\omega(G)-1$ if $\omega(G[N_3])\ge3$.

Therefore, we may suppose $\omega(G[N_3])\le2$. Then $\chi(G[\cup^8_{i=2}N_i])\le\chi(G[N_2\cup N_6])+\chi(G[N_3\cup N_7])+\chi(G[N_4\cup N_8])+\chi(G[N_5])\le \omega(G)-1+5\le2\omega(G)-1$.

\medskip

\noindent{\bf Case} 2. $\omega(G[N_2])\le3$.

Recall $\omega(G)\ge4$ by our assumption at the beginning of this section. If $\omega(G[N_2])\le2$, then $\chi(G[\cup^8_{i=2}N_i])\le\chi(G[N_2\cup N_6])+\chi(G[N_3\cup N_7])+\chi(G[N_4\cup N_8])+\chi(G[N_5])\le2+2+2+1=7\le2\omega(G)-1$. If $\omega(G[N_2])=3$ and $\omega(G[N_3])\le2$, then $\chi(G[\cup^8_{i=2}N_i])\le\chi(G[N_2\cup N_6])+\chi(G[N_3\cup N_7])+\chi(G[N_4\cup N_8])+\chi(G[N_5])\le3+2+2+1=8\le2\omega(G)$. If $\omega(G[N_2])=3, \omega(G[N_3])=3$ and $\omega(G[N_4])\le2$, then $G[N_2\cup N_3]$ is perfect by Claim \ref{t3-2} (3) and Observation \ref{hole}, which leads to $\chi(G[\cup^8_{i=2}N_i])\le\omega(G)+2+1\le2\omega(G)-1$. These follow that $\chi(G[\cup^8_{i=2}N_i])\le2\omega(G)$ if $\omega(G[N_2])\le3$.

This proves Claim \ref{t3-3}.
\qed

\medskip

Since $v_1$ is anticomplete to $N^{\ge2}(Q)$, it follows that $\chi(G[\{v_1\}\cup N^{\ge2}(Q)])=1$ as $N^{\ge2}(Q)$ is a stable set. Moreover, $\chi(G[N_1])\le\frac{1}{2}(\omega(G)^2-\omega(G))$ by Lemma \ref{claw}, and $\chi(G[\cup^8_{i=2}N_i])\le2\omega(G)$ by Claim \ref{t3-3}. We have that $\chi(G)\le \chi(G[N_1])+\chi(G[\cup^8_{i=2}N_i])+\chi(G[\{v_1\}\cup N^{\ge2}(Q)])\le \frac{1}{2}(\omega(G)^2-\omega(G))+2\omega(G)+1=\frac{1}{2}\omega^2(G)+\frac{3}{2}\omega(G)+1$.

This completes the proof of Therefore \ref{P3}.\qed

\medskip

\begin{thebibliography}{9999}
	
	\bibitem{BA18} A.P. Bharathi, S.A. Choudum, Colouring of ($P_3\cup P_2$)-free graphs, Graphs Comb., 34
	(2018) 97-107.
	
	\bibitem{BRS19} C. Brause, B. Randerath, I. Schiermeyer, E. Vumar, On the chromatic number of $2K_2$-free graphs, Disc. Applied Math., 253 (2019) 14-24.

    \bibitem{BM08} J.A. Bondy, U.S.R. Murty, Graph Theory, Springer, New York, 2008.
	
	\bibitem{CKS07} S. Choudum, T. Karthick, M. Shalu, Perfect coloring and linearly $\chi$-bound $P_6$-free graphs, J. of Graph Theory, 54 (2007) 293-306.
	\bibitem{CRSR06} M. Chudnovsky, N. Robertson, P. Seymour, R. Thomas, The strong perfect graph theorem, Annals of math., 164 (2006) 51-229.
	
	\bibitem{CS05} M. Chudnovsky, P. Seymour, The structure of claw-free graphs, in: Surveys in Combinatorics 2005, Cambridge Univ. Press, 2005, pp. 153-C171.
		
    \bibitem{CS10}	M. Chudnovsky, P. Seymour, Claw-free graphs VI. Colouring, J. of Combinatorial Theory, Ser. B, 100 (2010) 560-572.

    \bibitem{DXX22} W. Dong, B. Xu, Y. Xu, On the chromatic number of some $P_5$-free graphs, Disc. Math., 345 (2022) 113004.

    \bibitem{ER59}  P. Erd\"{o}s, Graph theory and probability, Canadian J. of Math., 11 (1959) 34-38.

    \bibitem{GM22}  M. Gei{\ss}er,  Colourings of $P_5$-free graphs, PhD thesis, 2022.

    \bibitem{Gy75} A. Gy\'{a}rf\'{a}s, On Ramsey covering-numbers, in: Colloquia Mathematic Societatis J\'{a}nos Bolyai 10, Infinite and Finite Sets, North-Holland/American Elsevier, New York, 1975, pp. 801-816.

    \bibitem{Gy87} A. Gy\'{a}rf\'{a}s, Problems from the world surrounding perfect graphs,  Applicationes Mathematicae, 19 (1987) 413-441.

    \bibitem{HCT18} C.T. Ho\'{a}ng, On the structure of (banner, odd hole)-free graphs, J. of Graph Theory, 89 (2018) 395-412.
	
	\bibitem{KKS22} T. Karthick, J. Kaufmann, V. Sivaraman, Coloring graph classes with no induced fork via perfect divisibility, The Electronic J. of Combinatorics, 29 (2022) \#P3.19.
	
	\bibitem{KIM95} J.H. Kim, The Ramsey number $R(3, t)$ has order of magnitude $O(\frac{t^2}{\log t})$, Random Structures and Algorithms, 7 (1995) 173-207.
	
	\bibitem{LW21} X. Liu, J. Schroeder, Z. Wang, X. Yu, Polynomial $\chi$-binding functions for $t$-broom-free graphs, arXiv preprint arXiv:2106.08871, 2021.
	
	\bibitem{PA22} A. Prashant, P. Francis, S.F. Raj, $\chi$-binding functions for some classes of $(P_3\cup P_2) $-free graphs, arXiv preprint arXiv:2203.06423, 2022.
	
	\bibitem{RB93} B. Randerath,  The Vizing bound for the chromatic number based on forbidden pairs, PhD thesis, RWTH Aachen, Shaker Verlag, 1993.
	
	\bibitem{RS04} B. Randerath, I. Schiermeyer, Vertex colouring and forbidden subgraphs-a survey, Graphs Comb., 20 (2004) 1-40.
	
	\bibitem{SR19} I. Schiermeyer, B. Randerath, Polynomial $\chi$-binding functions and forbidden induced subgraphs: a survey, Graphs Comb., 35 (2019) 1-31.
	
	\bibitem{SS20} A. Scott, P. Seymour, A survey of $\chi$-boundedness, J. of Graph Theory, 95 (2020) 473-504.
	
 	\bibitem{SS21} A. Scott, P. Seymour, S. Spirkl, Polynomial bounds for chromatic number. IV. A near-polynomial bound for excluding the five-vertex path, Combinatorica, to appear, arXiv:2110.00278.
	
	\bibitem{Su81} D.P. Sumner, Subtrees of a graph and chromatic number, in: The Theory and Applications of Graphs, John Wiley \& Sons, New York, 1981, pp. 557-576.
	
	\bibitem{Wa80} S. Wagon, A bound on the chromatic number of graphs without certain induced subgraphs, J. of Combinatorial Theory, Ser. B, 29 (1980) 345-346.
	
	\bibitem{WZ22} X. Wang, D. Zhang, The $\chi$-Boundedness of $(P_3\cup P_2)$-Free Graphs, J. of Mathematics, 2022.
	
	\bibitem{WX22} D. Wu, B. Xu, Perfect divisibility and coloring of some fork-free graphs, arXiv preprint arXiv:2302.06800, 2023.
	

	
\end{thebibliography}
\end{document}